\documentclass{amsart}
\usepackage{graphicx}
\usepackage[hyperindex,pdfmark]{hyperref}

\theoremstyle{definition}

\theoremstyle{remark}

\numberwithin{equation}{section}

\begin{document}

\title[Numerical study of the zeta derivative at zeros]{Numerical study of the derivative of the Riemann zeta function at zeros}

\author[G.A. Hiary]{Ghaith A. Hiary}
\thanks{Preparation of this material was partially supported by the National Science Foundation
under agreements No.  DMS-0757627 (FRG grant) and DMS-0635607.  Computations were carried out
at the Minnesota Supercomputing Institute.}
\address{School of Mathematics, University of Bristol, University Walk, Bristol, BS8 1TW.}
\email{hiaryg@gmail.com}

\author[A.M. Odlyzko]{Andrew M. Odlyzko}
\thanks{}
\address{School of Mathematics, University of Minnesota, 206 Church St. S.E., Minneapolis, MN, 55455.}
\curraddr{}
\email{odlyzko@umn.edu}

\subjclass[2000]{Primary, Secondary}
\maketitle

\begin{center}
{\em Dedicated to Professor Akio Fujii on his retirement.}
\end{center}

\begin{abstract}
The derivative of the Riemann zeta function
was computed numerically on several large sets of zeros at large heights.
Comparisons to known and conjectured
asymptotics are presented. 
\end{abstract}

\section{Introduction}

Throughout this paper, we assume the truth of the Riemann Hypothesis (RH),
and we let $\gamma_n >0$ denote the ordinate of the $n$-th non-trivial zero of $\zeta(s)$.  
Hejhal~\cite{He} assumed the RH and a weak consequence of Montgomery's~\cite{Mo} 
pair-correlation conjecture, 
namely that for some $\tau>0$, there is a constant $B$ such that
\begin{equation}
\limsup_{N\to \infty}\frac{1}{N}\left|\{n \,:\, N\le n\le 2N\,,\, (\gamma_{n+1}-\gamma_n)\log\gamma_n < c\}\right|\le B c^{\tau}\,, 
\end{equation}

\noindent
holds for all $c\in (0,1)$.  Under these assumptions, he
proved the following central limit theorem: for $\alpha <\beta$,  
\begin{equation} \label{eq:derclt1}
\lim_{N\to \infty} \frac{1}{N}\,\left|\left\{ n\,:\, N\le n\le 2N\,,\, \frac{\log \left|\displaystyle \frac{2\pi \zeta'(1/2+i\gamma_n)}{\log (\gamma_n/2\pi)}\right|}{\sqrt{\frac{1}{2} \log \log N}}\in(\alpha,\beta)\right\}\right| =
\frac{1}{\sqrt{2\pi}} \int_{\alpha}^{\beta} e^{-x^2/2}\,dx\,,
\end{equation}

\noindent
So under these assumptions, $\log|\zeta'(1/2+i\gamma_n)|$, 
suitably normalized, converges in distribution over fixed ranges to a  standard normal variable.
To obtain more precise information about the tails of the distribution,
we consider the moments

\begin{equation} \label{eq:dermo1}
J_{\lambda}(T) := \frac{1}{N(T)} \sum_{0<\gamma_n\le T} |\zeta'(1/2+i\gamma_n)|^{2\lambda}\,,
\end{equation}

\noindent
where $N(T):=\sum_{0<\gamma_n\le T} 1= \frac{T}{2\pi}\log \frac{T}{2\pi e}+ O(\log T)$
is the zero counting function. 
Notice that $J_{\lambda}(T)$ is defined for all $\lambda$ provided the zeros of $\zeta(s)$ 
are simple, as is widely believed.

Gonek~\cite{Go1}~\cite{Go2} carried out an extensive study of $J_{\lambda}(T)$.
He proved, under the assumption of the RH, that $J_1(T)\sim \frac{1}{12} (\log T)^3$ as
$T\to\infty$. It was suggested by Gonek~\cite{Go2}, and independently by Hejhal~\cite{He}, 
that $J_{\lambda}(T)$ is on the order of $(\log T)^{\lambda(\lambda+2)}$. 
 Ng~\cite{Ng} proved, under the RH, that $J_2(T)$ is  order of $(\log T)^8$, 
 which is in agreement with that suggestion. 

Hughes, Keating, and O'Connell~\cite{HKO}, applied the random matrix 
philosophy (e.g. see~\cite{KS}), which 
predicts that certain behaviors of $L$-functions are mimicked statistically by characteristic
polynomials of large matrices from the classical compact groups.  This led them to   
predict that for Re($\lambda$) $>-3/2$,

\begin{equation}\label{eq:hkc}
J_{\lambda}(T)\sim a(k)\frac{G^2(\lambda+2)}{G(2\lambda+3)} 
\left(\log \frac{T}{2\pi}\right)^{\lambda(\lambda+2)} ~~~~~~{\hbox{as}} ~~~~~~~~~~T \to \infty,
% \,,\qquad \Re(\lambda)>-3/2\,,
\end{equation}

\noindent
where $G(z)$ is the Barnes G-function, and $a(k)$ is an ``arithmetic factor.''
The conjecture (\ref{eq:hkc}) is consistent with previous theorems and conjectures.
%(as can be checked
% upon noting $a(-1)=6/\pi^2$ (by taking a limit), $a(1)=1$, $G(1)=1$, and $G^2(3)/G(5)=1/12$.

Recently, Conrey and Snaith~\cite{CS}, assuming the ratios conjecture, 
gave lower order terms in asymptotic expansions for $J_1(T)$ and $J_2(T)$.
They conjectured the existence of certain polynomials $P_{\lambda}(x)$,
for $2\lambda=2$ and $2\lambda=4$, such that

\begin{equation}
\sum_{0<\gamma_n< T}|\zeta'(1/2+i\gamma_n)|^{2\lambda}\sim\int_0^T P_{\lambda}\left(\log \frac{t}{2\pi}\right)\,dt\,,
\end{equation}

\noindent
The conjecture for the case $2\lambda=2$ was subsequently proved by Milinovich~\cite{Mi}, assuming the RH. 
It is expected that such polynomials exist for other integer values of $\lambda>0$ as well.
% Such polynomial can be expected to exist for other integer values of $\lambda>0$ as well.
%Since these results are precise they are very useful in numerical investigations.

%\begin{equation}\label{eq:csp12}
%\begin{split}
%P_1(x):=&\frac{1}{24\pi} x^4 +\frac{c_0}{3\pi}x^3+(\frac{c_0^2}{2\pi}-\frac{c_1}{\pi})x^2+\cdots\\
%P_2(x):=&\frac{1}{8640\zeta(2)}x^9+\frac{(2c_0\zeta(2)-\zeta'(2))}{480\zeta^2(2)}x^8\\
%        & +\frac{(7c_0^2 \zeta^2(2) - 2 c_1\zeta^2(2) -
%           8 c_0 \zeta(2) \zeta'(2) + 4 (\zeta'(2))^2 -
%           2  \zeta(2) \zeta''(2))}{120\zeta^3(2)}x^7+\cdots
%\end{split}
%\end{equation}
  
The purpose of this article is to study numerically various statistics  
 of the derivative of the zeta function at its zeros. In particular, we consider 
 the distribution of $\log|\zeta'(1/2+i\gamma_n)|$, 
  moments of $|\zeta'(1/2+i\gamma_n)|$, 
  and correlations among moments.
The goal is to obtain more detailed information about
 the derivative at zeros, and to enable comparison with various
 conjectured and known asymptotics. Our computations rely on large sets of zeros at large 
 heights that are described in detail in \cite{HO}.
 
We find that the empirical distribution of $\log|\zeta'(1/2+i\gamma_n)|$, 
  normalized to have mean zero and standard deviation one, agrees 
generally well with the limiting normal distribution proved by Hejhal, as shown in Figure~\ref{distder23}.
 But the empirical mean and standard deviation pre-normalization
 are noticeably different from predicted ones.
 Also, as shown in Figure~\ref{tailder23}, the frequency
 of  very small normalized values of $\log|\zeta'(1/2+i\gamma_n)|$ is higher than predicted by 
 a standard normal distribution, while the frequency of very large normalized values is lower than predicted. 
 Since these differences appear to 
 decrease steadily with height, however, they are probably not significant. 
 
To examine the tails of the distribution of $\log |\zeta'(1/2+i\gamma_n)|$, 
 we present data for the moments of $|\zeta'(1/2+i\gamma_n)|$ over short ranges:

\begin{equation}
J_{\lambda}(T,H):= \frac{1}{N(T+H)-N(T)} \sum_{T\le \gamma_n\le T+H} |\zeta'(1/2+i\gamma_n)|^{2\lambda} .
\end{equation}

\noindent
For large $\lambda$, the empirical values of $J_{\lambda}(T,H)$   
 deviate substantially from the values suggested by the
  leading term prediction (\ref{eq:hkc}).
 This is not surprising. Because 
 for $\lambda$ large relative to $T$, the contribution of lower order
 terms is likely to dominate, and so the leading term asymptotic on its own 
  may not suffice. 
  Furthermore, the said deviations decrease steadily with height  
  and they occur in a generally uniform way for roughly $2\lambda \le 6$,
 so they are consistent with the effect of ``lower order terms'' still 
 being felt even at such relatively large heights.

In the specific cases of the second and fourth moments of $|\zeta'(1/2+i\gamma_n)|$, 
 the conjectures of Conrey and Snaith~\cite{CS} supply lower order terms,  
 and the agreement with the data is much better, as shown in Table~\ref{csmomentder23}. 
 \footnote{It might be worth mentioning that 
 we attempted to calculate the coefficients of lower order terms 
 in the \cite{CS} conjectures by calculating $J_{\lambda}(T)$ 
for sufficiently many values of $T$, then solving the resulting system of equations. However,
this did not yield good approximations of the coefficients (even for small $\lambda$),
 which is not surprising, since the scale is logarithmic and the Conrey and Snaith expansion
is only asymptotic.} 
 
As $\lambda$ increases, the observed variability in the moments of $|\zeta'(1/2+i\gamma_n)|$ 
 is more extreme, but it is still
 significantly less than we previously encountered in the moments of $|\zeta(1/2+it)|$ (see \cite{HO}). 
 To illustrate, our computations of the twelfth moment of $|\zeta'(1/2+i\gamma_n)|$ over
 15 separate sets of $\approx 10^9$ zeros each (near the $10^{23}$-rd zero) 
show that the ratio of highest to lowest moment among the 15 twelfth moments thus obtained was 2.36. 
 In contrast, that ratio 
 for the twelfth moment of $|\zeta(1/2+it)|$ 
  was 16.34, which is significantly larger 
  (see \cite{HO}).
% observed  moment for $2\lambda = \, 6,\,8\,,10,$ and 12, 
% are, respectively, 1.03, 1.17, 1.54, and 2.36, which is a 
% (this can be deduced easily from Table~\ref{momentder23}).
% in contrast, the corresponding ratios for moment of $|\zeta(1/2+it)|$  
% computed over the intervals spanned by each these 15 samples,
% are significantly larger (they are 1.25, 2.44, 6.20, and 16.34, respectively; see~\cite{HO}). 
%(for $2\lambda= \,6,\,8,\,10,$ and 12, respectively).  

In general, the variability in statistical data for $|\zeta'(1/2+i\gamma_n)|$ 
is considerably less than the variability in  statistical data for $|\zeta(1/2+it)|$.
 It is not immediately clear why this should be so, considering, for instance, that 
 the central limit theorem for $\log|\zeta'(1/2+i\gamma_n)|$ is only conditional, while that for $\log|\zeta(1/2+it)|$ is not, 
 and both theorems scale by the same asymptotic variance.
% It might be interesting 
 %test such conjectured full asymptotics for higher $\lambda$ as well.

In the case of negative moments, our data is in agreement with Gonek's conjecture (\cite{Go1}) 
   $J_{-1}(T) \sim \frac{6}{\pi^2}(\log T/(2\pi))^{-1}$ as $T\to\infty$. 
  But starting at $2\lambda=-3$, 
 and as $\lambda$ decreases, the empirical behavior of
   negative moments becomes rapidly more erratic.
 For example, using the same 15 zero sets near the $10^{23}$-rd zero mentioned previously, 
  the ratio of highest to lowest negative moment among them 
 gets very large as $\lambda$ decreases; we obtain: 1.03, 8.45, 178.49, and 17240.99, for $2\lambda= -2, -3,-4,$
 and $-6$, respectively (this can be deduced easily from Table~\ref{Nmomentder23}). 
  Notice that the point $2\lambda=-3$ is special because it is where the leading term prediction
 (\ref{eq:hkc}) first breaks down 
 due to a pole of order 1 in the ratio of Barnes G-functions.

 Extreme values of negative
 moments are caused by very few zeros. When $2\lambda=-3$, for instance,
 the largest observed moment 
 among our 15 sets is $0.178047$. 
 About 87\% of this value is contributed by 4 zeros where 
 $|\zeta'(1/2+i\gamma_n)|$ is small and equal to 
  0.002439, 0.002453, 0.004388, and 0.004365.\footnote{We checked such small values of $|\zeta'(1/2+i\gamma)|$ 
 by computing them in two ways, using the Odlyzko-Sch\"onhage algorithm, 
 and using the straightforward Riemann-Siegel formula;  
 the results from the two methods agreed to within $\pm 10^{-6}$}
 Such small values of $|\zeta'(1/2+i\gamma)|$
  typically occur at pairs of consecutive zeros that are close to each other.
For example,
 the values 0.002439 and 0.002453 occur 
 at the following two consecutive zero ordinates:

\begin{equation}
\begin{split}
& 1.30664344087942265202071895041619 \times 10^{22}\,,\\
& 1.30664344087942265202071898265199 \times 10^{22}\,.
\end{split}
\end{equation}

\noindent
The above pair of zeros is separated by 0.00032, which is about 1/400 times
the average spacing of zeros at that height (which is $\approx 0.128$). 

%a5 G23 outder628
%10^(22)*1.306643440879414474 + 81780207 + 0.1895041619
%der 0.002439
%check y
%10^(22)*1.306643440879414474 + 81780207 + 0.1898265199
%der -0.002453
%check y

%a5 G23 outder76
%10^(22)*1.306643440879414474 +  9907257 + 0.3240763646
%der 0.004388
%check y
%10^(22)*1.306643440879414474 +  9907257 + 0.3244584141
%der -0.004365
%check y

To investigate possible correlations among  values of $|\zeta'(1/2+i\gamma_n)|^{2\lambda}$,
we studied numerically the (shifted moment) function: 
 
\begin{equation}
S_{\lambda}(T,H,m):= \sum_{T\le \gamma_n\le T+H} |\zeta'(1/2+i\gamma_n)\zeta'(1/2+i\gamma_{n+m})|^{2\lambda}\,.
\end{equation}

\noindent
 We plotted $S_{\lambda}(T,H,m)$, for several choices of  $\lambda$, $T$, and $H$, 
 and as $m$ varies. The resulting plots indicate there are
 long-range correlations among the values of the derivative at zeros. 
Unexpectedly, the tail of $S_2(T,H,m)$ (Figure~\ref{corrder23}; right plot) strongly resembles  
 the tail for the shifted fourth moment of $|\zeta(1/2+it)|$ (Figure~4 in \cite{HO}).

To better understand these correlations, we considered the ``spectrum'' of $\log|\zeta'(1/2+i\gamma_n)|$; 
 see (\ref{eq:specf}) for a definition. 
 A plot of the spectrum reveals sharp spikes, shown in Figure~\ref{fftder23}.
 These spikes can be explained heuristically by applying techniques already used
by Fujii \cite{Fu,Fu2} and Gonek \cite{Go1} to estimate sums involving $\zeta'(1/2+i\gamma_n)$.

\section{Numerical results}

 Conjecture (\ref{eq:derclt1}) suggests the mean and standard deviation
  of $\log |\zeta'(1/2+i\gamma)|$ for 
 zeros from near $T=1.3066434\times 10^{22}$ (i.e. near the $10^{23}$-rd zero) 
 are about 2.0 and 1.4, respectively. This is far from the empirical mean and standard deviations
 listed in Table~\ref{sum0}, which are 3.4907 and 1.0977.
\footnote{The mean and standard deviations listed in Table~\ref{sum0} change very little across different
 zero sets near the same height. For example, using a different set of $10^8$ zeros near 
 the $10^{23}$-rd zero, the empirical mean is 3.4907 and the empirical
 standard deviation is 1.0978, which are very close the numbers listed in Table~\ref{sum0}.
We note that the empirical mean and standard deviation 
  are closer to the values suggested by the central limit theorem for characteristic polynomials of
 unitary matrices (see \cite{HKO}), which are 3.47 and 1.12.} 
 Since these quantities grow very slowly (like $\log \log T$),
 these differences are probably not significant. 
%using 25M zeros: mean = 3.4908, sd = 1.0978

%rmt mean:
%3.0924
%3.3227
%3.4707

%hejhal mean:
%1.677318
%1.907622
%2.055706

%rmt sd:
%1.038655
%1.092683
%1.126054

%hejhal sd:
%1.325744
%1.368484
%1.395275

\begin{table}[ht]
\footnotesize
\caption{\footnotesize Summary statistics for  $\log|\zeta'(1/2+i\gamma_n)|$ using sets of $10^7$ zeros from different heights 
 The column "Zero" lists the zero number near which the set is located. 
 SD stands for standard deviation. \label{sum0}}
\begin{tabular}{|c|c|c|c|c|}
\hline
Zero &Min & Max & Mean & SD \\
\hline
$10^{16}$ & -3.7371 & 7.3920  & 3.1211 & 1.0135 \\
$10^{20}$ & -3.2181 & 8.0085  & 3.3458 & 1.0653 \\ 
$10^{23}$ & -2.9602 & 8.2836  & 3.4907 & 1.0977 \\  
\hline
\end{tabular}
\end{table}

\begin{figure}[ht]
\footnotesize
\caption{\footnotesize Empirical density of $\log |\zeta'(1/2+i\gamma_n)|$, after being normalized to have mean 
0 and standard deviation 1, using $10^7$ values of $\log |\zeta'(1/2+i\gamma_n)|$ from near 
the $10^{23}$-rd zero (the bin size is $0.0512$). The density of a standard normal 
variable (continuous line) is drawn to facilitate comparison.  The right plot shows the
diference, the empirical density minus the normal.}\label{distder23}
%min -5.88, max 4.37, bin size 0.0512}
\includegraphics[scale=.27]{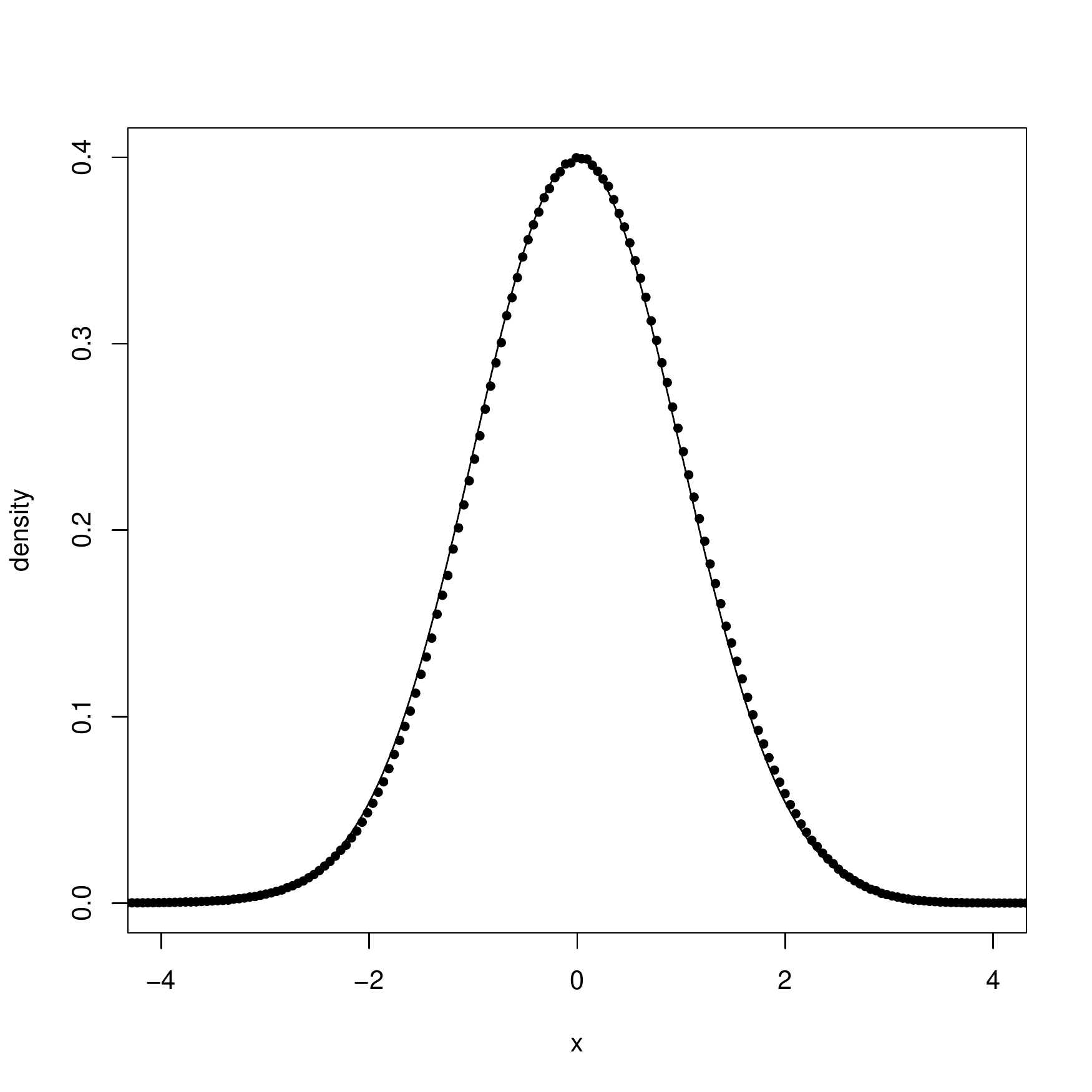}
\includegraphics[scale=.27]{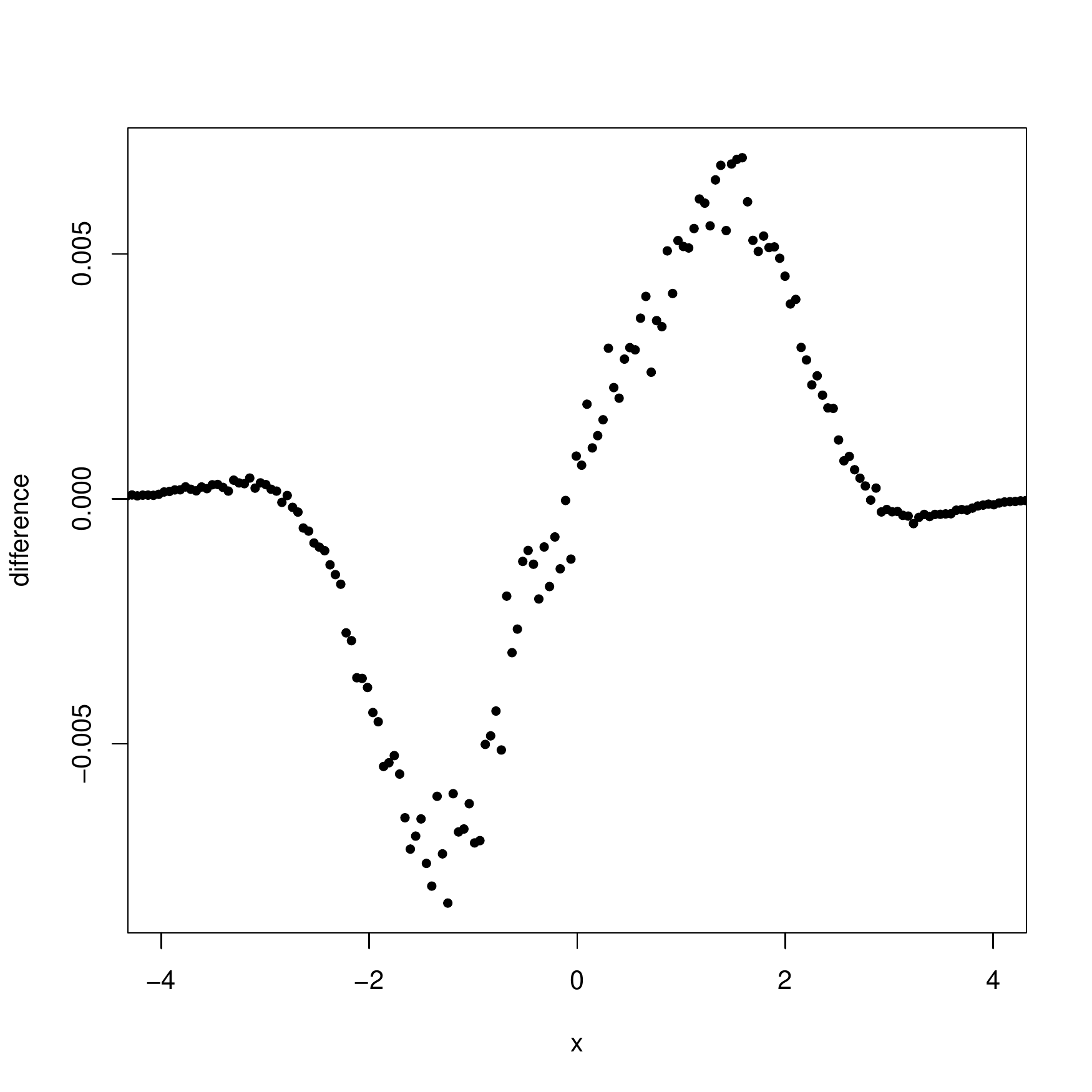}
\end{figure}

We normalize the sequence $\{\log |\zeta'(1/2+i\gamma_n)|\,:\, N \le n \le N + 10^7\}$, where $N\approx 10^{23}$,
to have mean zero and variance one. 
The distribution of the normalized sequence is illustrated in Figure~\ref{distder23}, which contains
 two plots, one of the empirical density function, and another of the difference between the empirical
 density and the predicted (standard Gaussian) density $\frac{1}{\sqrt{2\pi}}e^{-x^2/2}$.
 The fit in the first plot is visibly good, but there is a slight shift to the right about the center. 
 This shift is made more visible in 
 the second plot, which shows that the empirical density is generally larger than expected for $x>0$, 
 and is smaller than expected for $x<0$. 

\begin{figure}[ht]
\footnotesize
\caption{\footnotesize Distribution at the tails using $1.5\times 10^{10}$ zeros near the $10^{23}$-rd zero (bin size is 0.01).}\label{tailder23}
%left tail: min -8.66, max -3.18, bin size 0.01. right tail: min 2.98, max 4.89, bin size 0.01}
\includegraphics[scale=.27]{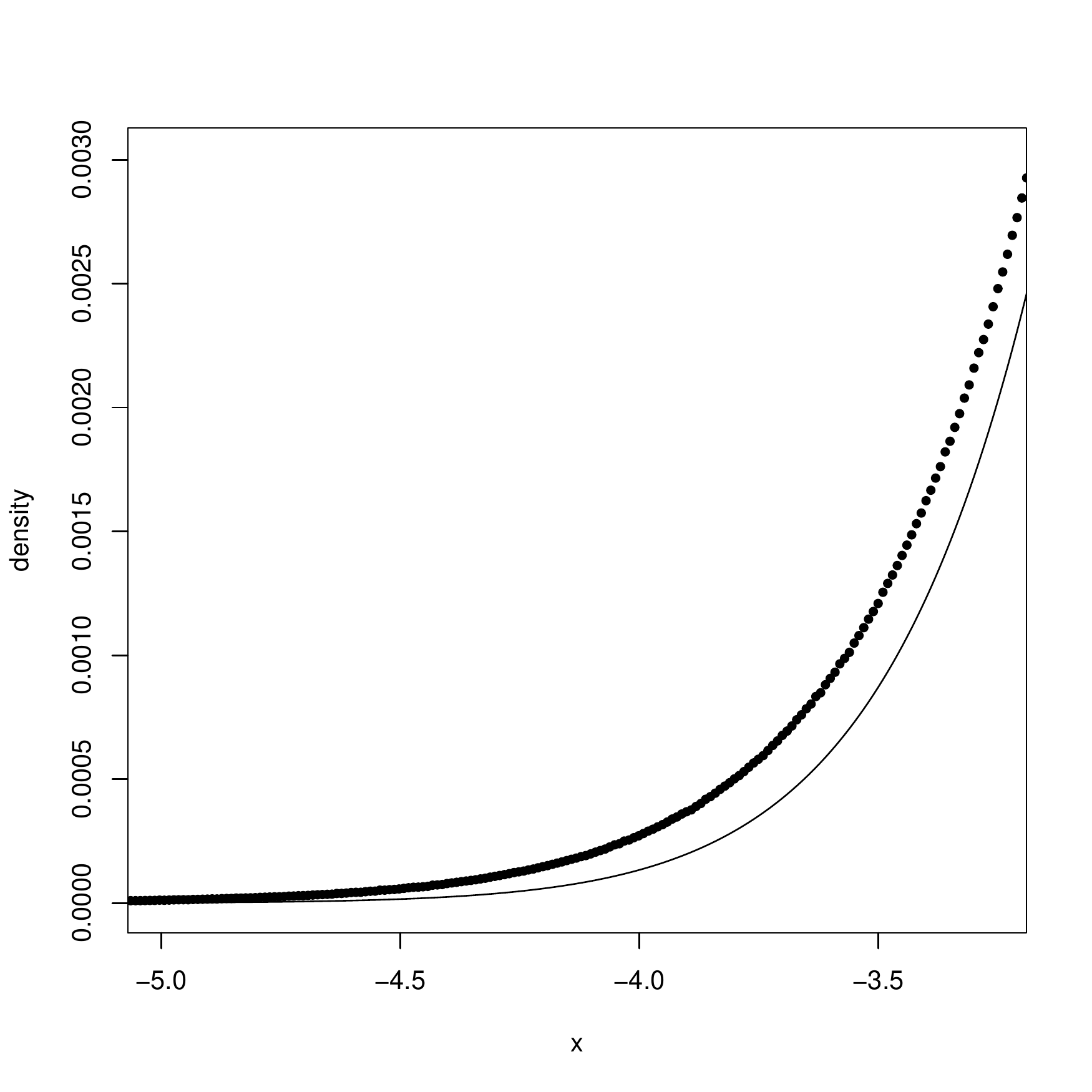}
\includegraphics[scale=.27]{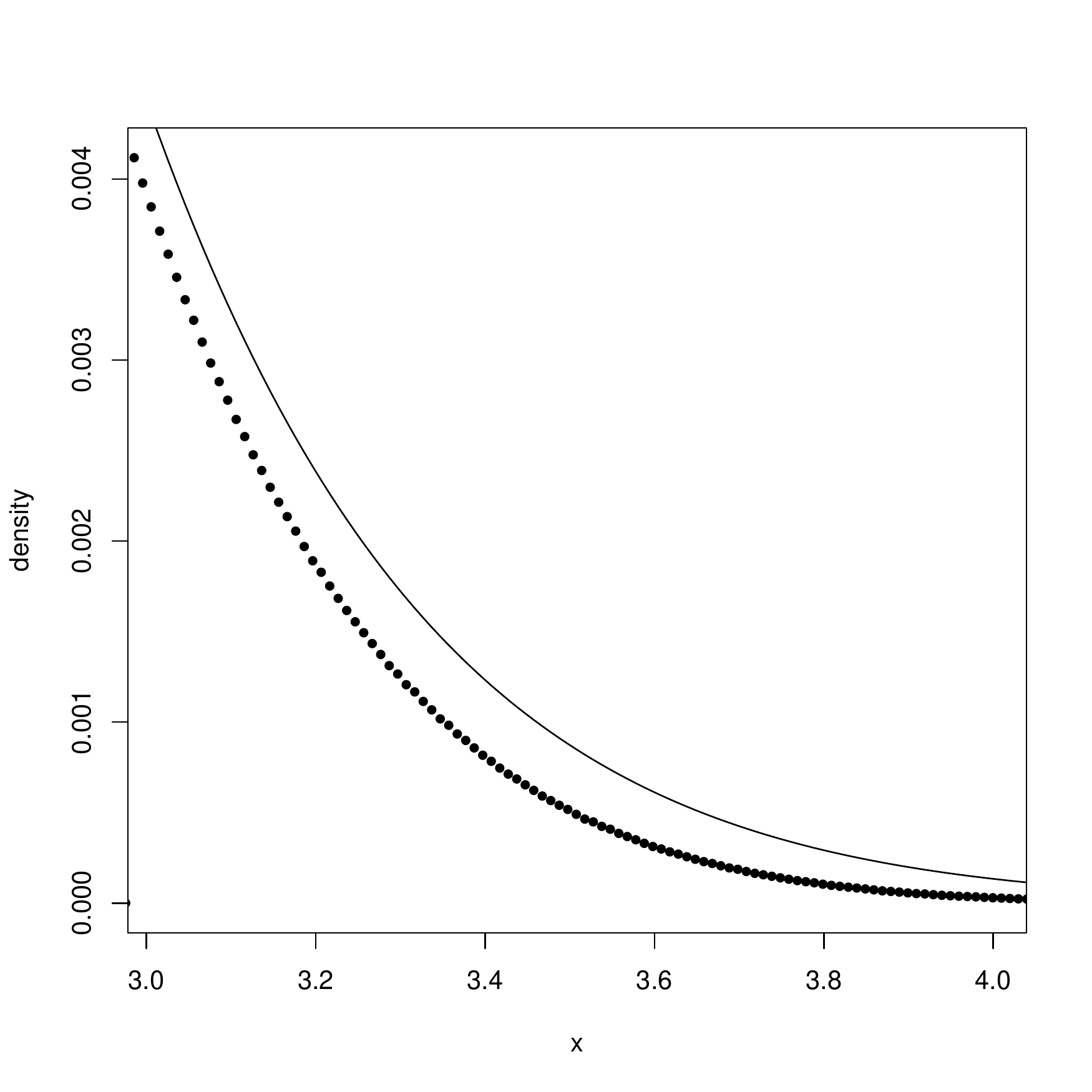}
\end{figure}

Near the tails, however, the situation is reversed. Figure~\ref{tailder23} 
 shows there is a deficiency in the occurrence of very large values of $|\zeta'(1/2+i\gamma_n)|$, and an abundance in the occurrence of very small values. 
For instance, conjecture (\ref{eq:derclt1}) suggests that about 
0.1462\% of the values of $|\zeta'(1/2+i\gamma_n)|$ near the $10^{23}$-rd zero  
 should satisfy $|\zeta'(1/2+i\gamma_n)| > 860$, 
 which is noticeably larger than the observed 0.1056\%.
 The conjecture also suggests about 0.0736\% of the values should satisfy $|\zeta'(1/2+i\gamma_n)|<1$, which 
 is smaller than the observed 0.1051\%.

 We remark the behavior near the tails becomes more consistent with
 expectation as height increase . For example, only 0.0025\% of the time
do we have $\log|\zeta'(1/2+i\gamma_n)|>3.2$ near the $10^{16}$-th zero,
which is far from the expected 0.068\%, but the percentage increases to 0.040\% near the $10^{23}$-rd zero.

 %max23: 7057.970000 6907.027561 6658.674692 6636.362297 6399.952231
%min23: 0.002439 0.002453 0.002719 0.002737 0.003094

%normal( > 3.2)  0.000687   
%15B   ( > 3.2)  0.000405
%1B    ( > 3.2)  0.000405 
%O20   ( > 3.2)  0.000362    
%Z16   ( > 3.2)  2.556778e-05  

%normal( < -3.2) 0.000687
%15B   ( < -3.2) 0.000997
%1B    ( < -3.2) 0.000996
%O20   ( < -3.2) 0.001030
%Z16   ( < -3.2) 0.001095 

%15B(|zeta'| > 860)  0.00105654
%15B( |zeta'| < 1)   0.00105129
%normal( < -3.180012) 0.00073634
%normal( > 2.975523)  0.00146245 

%xfit0<-seq(min(data23),max(data23),length=200); h0<-hist(data23,breaks=c(min(data23)-0.01,xfit0),plot=FALSE); plot(xfit0,h0$density, xlim=c(-4,4), type='p',pch=20,cex=1, xlab='x',ylab='density'); lines(xfit0,exp(-xfit0^2/2)/sqrt(2*pi)); min(data23); max(data23);(max(data23)-min(data23))/200

%USING EMPIRICAL MEAN AND SD BELOW
%--------------------------------------------
%normal: 0.0668072, 0.0013499, 3.39767e-6
%10M( < -1.5) 0.0661842, ( < -3.0) 0.001769, ( < -4.5) 2.0000e-5
%25M( < -1.5) 0.0662288, ( < -3.0) 0.001756, ( < -4.5) 1.773599e-5

%normal: 0.0668072, 0.0013499, 3.39767e-6
%10M( > 1.5) 0.0674187, ( > 3.0) 0.000944, ( > 4.5) 0
%25M( > 1.5) 0.0674019, ( > 3.0) 0.000949, ( > 4.5) 3.93259e-8

%normal( < -0.75) 0.22662
%10M23 ( < -0.75) 0.22561
%10M20 ( < -0.75) 0.22546 
%10M16 ( < -0.75) 0.22555 

%normal( > 0.75) 0.226627
%10M23 ( > 0.75) 0.226949
%10M20 ( > 0.75) 0.227223
%10M16 ( > 0.75) 0.227524 
%-------------------------------------------
%normal( > 3.2)  0.000687   
%15B   ( > 3.2)  0.000405
%1B    ( > 3.2)  0.000405 
%O20   ( > 3.2)  0.000362    
%Z16   ( > 3.2)  2.556778e-05  

%normal( < -3.2) 0.000687
%15B   ( < -3.2) 0.000997
%1B    ( < -3.2) 0.000996
%O20   ( < -3.2) 0.001030
%Z16   ( < -3.2) 0.001095 
%--------------------------------------------

%normal(  < -3.1407115) 0.0008426 
%O20( |zeta'| < 1)    0.0012195

%normal( < -3.079526)   0.0010366   
%ZZ16( |zeta'| < 1)    0.0015278

%normal( > 3.202039)    0.0006822  
%O20(|zeta'| > 860)   0.0003622

%normal( > 3.587402)    0.0001669  
%ZZ16(|zeta'| > 860)  2.5567783e-05
%-------------------------

For another  measure of the quality of the fit to the standard Gaussian in Figure~\ref{distder23}, 
we compare moments of both distributions. Table~\ref{gaussder23} shows the first few moments (the even moments in particular)
  agree reasonably well. Notice the odd moments tend to be negative, which is likely due to the aforementioned bias in the
 frequency of very small and very large values. 
  
\begin{table}[ht]
\footnotesize
\caption{\footnotesize Moments of $\log|\zeta'(1/2+i\gamma_n)|$, after being normalized to have mean zero and variance one,
 calculated using $10^7$ zeros from near the $10^{23}$-rd zero. The third column is the moment of a standard Gaussian.}\label{gaussder23}
\begin{tabular}{c|cc}
Moment & Derivatives & Gaussian \\
\hline
3rd& -0.02728 & 0\\
4th&  3.01364 & 3\\
5th& -0.49120& 0\\
6th&  15.3053& 15\\
7th& -7.43073& 0 \\
8th&  112.013& 105 \\
9th& -118.588& 0\\
10th & 1116.64 & 945
\end{tabular}
\end{table}

%-------------------------------
%max23: 7057.970000 6907.027561 6658.674692 6636.362297 6399.952231
%min23: 0.002439 0.002453 0.002719 0.002737 0.003094
%--------------------------

To better understand the tails of the distribution of $\log |\zeta'(1/2+i\gamma_n)|$, 
we consider the moments $J_\lambda(T)$ defined in (\ref{eq:dermo1}). 
 Since we are interested in the asymptotic behavior of $J_{\lambda}(T)$, 
we compare against the leading term prediction (\ref{eq:hkc}). 
We calculated ratios of the form  

\begin{equation}\label{eq:appder1}
\frac{\frac{1}{|B|}\sum_{\gamma\in B} |\zeta'(1/2+i\gamma)|^{2\lambda}}{a(\lambda) \frac{G^2(\lambda+2)}{G(2\lambda+3)} \left(\log \frac{T}{2\pi}\right)^{\lambda(\lambda+2)}}\,,
\end{equation}

\noindent
where $B$ is a block of consecutive zeros, $|B|$ denotes the number of zeros in $B$,
 and $T$ is the height where block $B$ lies.
If $T$ is large enough, one expects the value of (\ref{eq:appder1}) 
 to approach 1 as the block size $|B|$ increases.  
 Table~\ref{momentder23}, which uses blocks of size $|B|\approx 10^9$ (except for the first set, which uses the first $10^8$ zeta zeros), 
 shows that the empirical moments are significantly larger than the corresponding predictions, even for low moments. 
 For example, the empirical second moments ($2\lambda = 2$) near the $10^{23}$-rd zero  
  are generally off from expectation by about 9.6\%.   
 
 Nevertheless, the ratios (\ref{eq:appder1}) appear to decrease towards the expected 1 as the height increases,  
 and there is relatively little variation in the moment data for sets from near the same height when $2\lambda\le 6$.
 Both of these observations are consistent with the ``lower order terms'' still contributing significantly.

\begin{table}[ht]
\footnotesize
\caption{\footnotesize Ratio (\ref{eq:appder1}) calculated with $|B|\approx 10^9$, 
except for the first set, which uses the first $10^8$ zeros. The column ``Zero'' lists
 the approximate zero number near which block $B$ is located.} \label{momentder23}
\begin{tabular}{|l||l|l|l|l|l|l|}
\hline
Zero & $2\lambda=2$ & $2\lambda=4$ & $2\lambda=6$ & $2\lambda=8$ & $2\lambda=10$ & $2\lambda=12$ \\
\hline
$10^8$ & 1.1247 & 3.1579 & 91.856 & 78341 & $4.1016\times 10^9$ & $2.3478\times 10^{16}$ \\ 
$10^{16}$ & 1.1424 & 2.2087 & 17.686 & 1266.9 & $1.5057\times 10^6$ & $4.9628\times 10^{10}$\\ 
$10^{20}$ & 1.1123 & 1.9102 & 10.943 & 422.72 & $1.9904\times 10^5$ & $1.8362\times 10^9$ \\
$10^{23}$ & 1.0964 & 1.7645 & 8.4406 & 233.63 & $6.4583\times 10^4$ & $2.7127\times 10^8$  \\
- & 1.0964 & 1.7603 & 8.1602 & 199.18 & $4.1647\times 10^4$ & $1.1369\times 10^8$ \\
- & 1.0964 & 1.7598 & 8.1879 & 202.40 & $4.3355\times 10^4$ & $1.2325\times 10^8$ \\ 
- & 1.0964 & 1.7629 & 8.3221 & 217.58 & $5.2539\times 10^4$ & $1.7809\times 10^8$ \\
- & 1.0964 & 1.7630 & 8.3861 & 228.51 & $6.2549\times 10^4$ & $2.6614\times 10^8$ \\
- & 1.0964 & 1.7600 & 8.2022 & 206.36 & $4.6423\times 10^4$ & $1.4200\times 10^8$ \\
- & 1.0965 & 1.7642 & 8.3321 & 218.38 & $5.3663\times 10^4$ & $1.8923\times 10^8$ \\
- & 1.0965 & 1.7612 & 8.1862 & 201.43 & $4.3256\times 10^4$ & $1.2547\times 10^8$\\
- & 1.0963 & 1.7590 & 8.2176 & 209.97 & $4.8853\times 10^4$ & $1.5596\times 10^8$ \\
- & 1.0964 & 1.7654 & 8.3856 & 217.09 & $4.8781\times 10^4$ & $1.4148\times 10^8$ \\
- & 1.0963 & 1.7616 & 8.3009 & 218.92 & $5.4691\times 10^4$ & $1.9491\times 10^8$ \\
- & 1.0964 & 1.7585 & 8.1576 & 204.55 & $4.6872\times 10^4$ & $1.5134\times 10^8$ \\
- & 1.0965 & 1.7615 & 8.2380 & 209.26 & $4.7946\times 10^4$ & $1.5078\times 10^8$ \\
- & 1.0963 & 1.7586 & 8.1764 & 203.00 & $4.4241\times 10^4$ & $1.2904\times 10^8$ \\
- & 1.0964 & 1.7603 & 8.2037 & 208.39 & $4.9019\times 10^4$ & $1.6822\times 10^8$ \\
\hline
\end{tabular}
\end{table}

The full moment prediction of \cite{CS}, which takes lower order terms into account,
might lead one to expect that for $2\lambda=2$, $2\lambda=4$, 
  as $T\to \infty$, and for blocks $B$ not too small compared to $T$, 

\begin{equation}\label{eq:dercs}
\sum_{\gamma\in B} |\zeta'(1/2+i\gamma)|^{2\lambda}\sim\int_{B} P_{\lambda}(\log (t/2\pi))\,dt\,,
\end{equation}

\noindent
 where $P_{\lambda}(x)$ is as given in \cite{CS}, 
and $\int_{B}$ is short for integrating over the interval spanned by the block $B$.
To test this, we calculated ratios of the form

\begin{equation}\label{eq:appder2}
\frac{\sum_{\gamma\in B} |\zeta'(1/2+i\gamma)|^{2\lambda}}{\int_{B} P_{\lambda}(\log (t/2\pi))\,dt}\,.
\end{equation}

\noindent
As the block size increases, we expect (\ref{eq:appder2}) to be significantly closer to 1 than (\ref{eq:appder1}) since
it relies on a more accurate prediction.
 This is indeed what Table~\ref{csmomentder23} illustrates, where we see the fit to moment data is much better
 than we found in Table~\ref{momentder23}. 
\footnote{Notice if $T$ is large compared to the length of the interval spanned by block $B$,
 the denominator in ratio (\ref{eq:appder2}) is largely a function of $T$ multiplied 
by the length of the interval spanned by $B$.} (We point out that in the case $2\lambda=4$  
 only the first three terms in the full moment conjecture were used, because these were the only terms
 provided explicitly in \cite{CS}. It is likely the fit to the data 
 will be even better if the missing terms are included.)

\begin{table}[ht]
\footnotesize
\caption{\footnotesize Ratio (\ref{eq:appder2}) calculated with $|B|\approx 10^9$, 
except for the first set, which uses the first $10^8$ zeros. The column ``Zero'' lists
 the approximate zero number near which block $B$ is located.}\label{csmomentder23}
\begin{tabular}{|l||l|l|}
\hline
Zero & $2\lambda=2$ & $2\lambda=4$ \\
\hline
$10^8$ &  1.0000 & 1.0924 \\
$10^{16}$ & 1.0000 & 1.0144 \\
$10^{20}$ & 1.0000 & 1.0087 \\
$10^{23}$ & 1.0000 & 1.0074 \\
~~~~`` & 1.0000 & 1.0050 \\
~~~~`` & 0.9999 & 1.0047 \\
~~~~`` & 1.0000 & 1.0064 \\
~~~~`` & 0.9999 & 1.0065 \\
~~~~`` & 0.9999 & 1.0048 \\
~~~~`` & 1.0000 & 1.0072 \\
~~~~`` & 1.0000 & 1.0055 \\
~~~~`` & 0.9998 & 1.0042 \\ 
~~~~`` & 1.0000 & 1.0079 \\
~~~~`` & 0.9999 & 1.0057 \\
~~~~`` & 0.9999 & 1.0039 \\
~~~~`` & 1.0000 & 1.0057 \\
~~~~`` & 0.9999 & 1.0040 \\
~~~~`` & 1.0000 & 1.0049 \\
\hline
\end{tabular}
\end{table}

We remark that the five largest values of $|\zeta'(1/2+i\gamma_n)|$ in our data set 
are $\approx$ 7057, 6907, 6658, 6636, and 6399.
 The cumulative contribution of these large values to the $2\lambda$-th moment, as a percentage of the overall $2\lambda$-th moment, is listed
 in Table~\ref{accumder23} for several $\lambda$. 

\begin{table}[ht]
\footnotesize
\caption{\footnotesize Cumulative contribution percentage of the 5 largest values of $|\zeta'(1/2+i\gamma_n)|$ to 
the empirical $2\lambda$-th moment
for $1.5\times 10^{10}$ zeros near the $10^{23}$-rd zero.}\label{accumder23}
\begin{tabular}{|l|l|l|}
\hline
$2\lambda=8$ & $2\lambda=10$ & $2\lambda=12$ \\
\hline
0.50 & 1.84  & 4.51 \\
0.92 & 3.32  & 7.99 \\
1.24 & 4.35  & 10.2 \\
1.54 & 5.35  & 12.3 \\
1.77 & 6.04  & 13.7 \\
\hline
\end{tabular}
\end{table}

In the case of negative moments, the conjecture $J_{-1}(T) \sim \frac{6}{\pi^2}(\log T/(2\pi))^{-1}$ as $T\to\infty$, due to Gonek~\cite{Go2},
 suggests the negative second moment
 should be $\approx 0.01808$ near zero number $10^{16}$, 
  $\approx 0.01436$ near zero number $10^{20}$, and $\approx 0.01238$ near zero number $10^{23}$.
 These predictions are in good agreement with the values listed in Table~\ref{Nmomentder23}. 
%is around 0.03864579 0.01808097 0.01436158 0.01238482

For $2\lambda \le -3$, the behavior 
 is much less predictable because, empirically, their sizes are determined by 
 a few zeros where $|\zeta'(1/2+i\gamma_n)|$ is small. 
 In fact, the particularly large fluctuations in the size of the negative
 sixth moment ($2\lambda = -6$), near the $10^{23}$-rd zero in Table~\ref{Nmomentder23}, are essentially due to 8 zeros (out of $1.5\times 10^{10}$)
 where $|\zeta'(1/2+i\gamma_n)|$ is equal to 0.002439, 0.002453, 0.002719, 0.002737, 0.003094, 0.003108, 0.004365, and 0.004388. 

\begin{table}[ht]
\footnotesize
\caption{\footnotesize Ratio (\ref{eq:appder2}) calculated with $|B|\approx 10^9$, 
except for the first set, which uses the first $10^8$ zeros. The column ``Zero'' lists
 the approximate zero number near which block $B$ is located.} \label{Nmomentder23}
\begin{tabular}{|l||l|l|l|l|}
\hline
Zero & $2\lambda=-2$ & $2\lambda=-3$ & $2\lambda=-4$ & $2\lambda=-6$\\
\hline
$10^8$ & 0.041129 & 0.059025 & 1.04212 & 2935.6 \\
$10^{16}$ & 0.018057 & 0.030660 & 0.55588 & 1488.1 \\
$10^{20}$ & 0.014341 & 0.028403 & 0.73586 & 2873.2 \\
$10^{23}$ & 0.012347 & 0.022040 &  0.41441 &    1106.5 \\
~~~~`` & 0.012365 & 0.022605 &  0.43869 &    1314.6 \\
~~~~`` & 0.012462 & 0.037677 &  2.76255 &    63336  \\
~~~~`` & 0.012321 & 0.021618 &  0.42275 &    1431.0  \\
~~~~`` & 0.012776 & 0.178047 &  59.6610 & 9288238  \\
~~~~`` & 0.012326 & 0.021062 &  0.33853 &     665.29 \\
~~~~`` & 0.012515 & 0.052929 &  7.46570 &  412318 \\
~~~~`` & 0.012334 & 0.022429 &  0.56305 &    4157.4 \\
~~~~`` & 0.012376 & 0.025800 &  0.81652 &    5414.6 \\
~~~~`` & 0.012541 & 0.089163 & 21.5695 & 2174342 \\
~~~~`` & 0.012411 & 0.039415 &  4.32860 &  185114 \\
~~~~`` & 0.012329 & 0.022729 &  0.55154 &    2723.6 \\
~~~~`` & 0.012386 & 0.027487 &  1.08706 &   11563 \\
~~~~`` & 0.012605 & 0.117993 &  35.4067 & 4686740 \\
~~~~`` & 0.012334 & 0.021217 &  0.33424 &     538.73 \\
\hline
\end{tabular}
\end{table}

%\begin{figure}[ht]
%\footnotesize
%\caption{\footnotesize Shifted 2nd moments, $10^7$ zeros near 10e23.}\label{corrder231}
%\includegraphics[scale=.3]{moment23_shift2_1}
%\includegraphics[scale=.3]{moment23_shift2_2}
%\includegraphics[scale=.3]{moment23_shift2_3}
%\end{figure}

%> pdf('moment23_shift2_1.pdf')
%> xfit<-seq(0,25,by=1); plot(xfit,cv23_1$acf[xfit+1]/cv23_1$acf[1], type='l', lwd=1.5, xlab='m',ylab='c_m/c_0'); text(18,0.7,"c_0 = 10792.7")
%> dev.off()
%> pdf('moment23_shift2_2.pdf')
%> xfit<-seq(3,325,by=1); plot(xfit,cv23_1$acf[xfit+1]/cv23_1$acf[1], type='l', lwd=1.5, xlab='m',ylab='c_m/c_0',xaxt='n'); axis(side=1,at=c(3,seq(50,300,length=6)))
%> dev.off()
%> pdf('moment23_shift2_3.pdf')
%> xfit<-seq(800,999,by=1); plot(xfit,cv23_1$acf[xfit+1]/cv23_1$acf[1], type='l', lwd=1.5, xlab='m',ylab='c_m/c_0',xaxt='n'); axis(side=1,at=c(800,850,900,950,1000))
%> dev.off()

Starting with the investigations of \cite{Od2}, several long-range correlations
have been found experimentally in zeta function statistics.  Such correlations are
not present in random matrices, but do appear in some dynamical systems that
for certain ranges are modeled by random matrices.  So far all the zeta function
correlations of this nature have been explained (at least numerically and heuristically)
by relating them to known properties of the zeta function, such as explicit
formulas that relate primes to zeros.  A natural question is whether such
correlations arise among values of $\zeta'(1/2+i\gamma_n)$.

In order to detect correlations among values of $|\zeta'(1/2+i\gamma_n)|$, 
consider

 \begin{equation}
 S_{2}(T,H,m):= \sum_{T\le \gamma_n\le T+H} |\zeta'(1/2+i\gamma_n)\zeta'(1/2+i\gamma_{n+m})|^4\,.
 \end{equation}

\begin{figure}[ht]
\footnotesize
\caption{\footnotesize Plots of $S_2(T,H,m)/S_2(T,H,0)$ using $10^7$ zeros near the $10^{16}$-th (left plot) and the $10^{23}$-rd zero (right plot)}\label{corrder23}
%up to 222 and 325 resp.
\includegraphics[scale=.3]{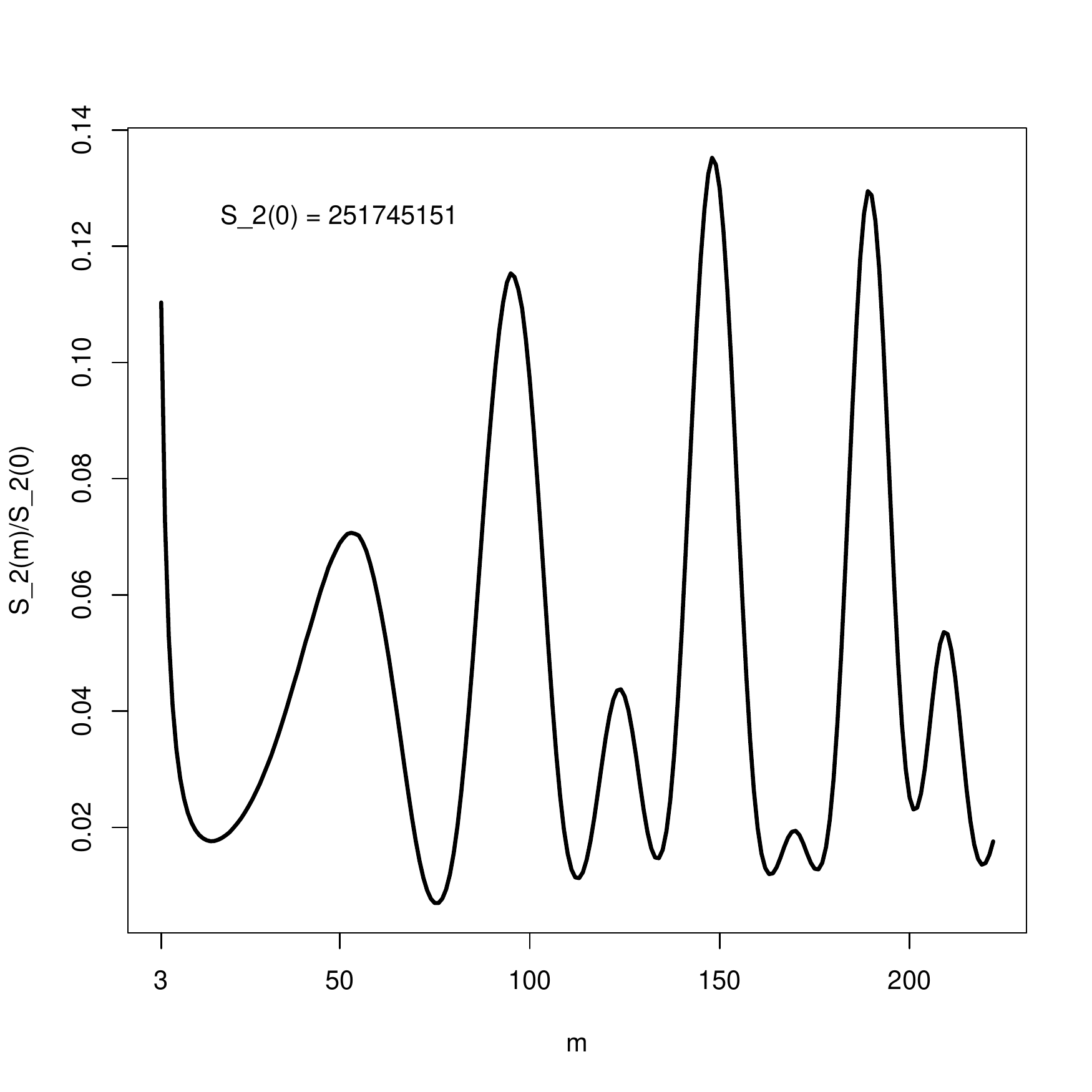}
\includegraphics[scale=.3]{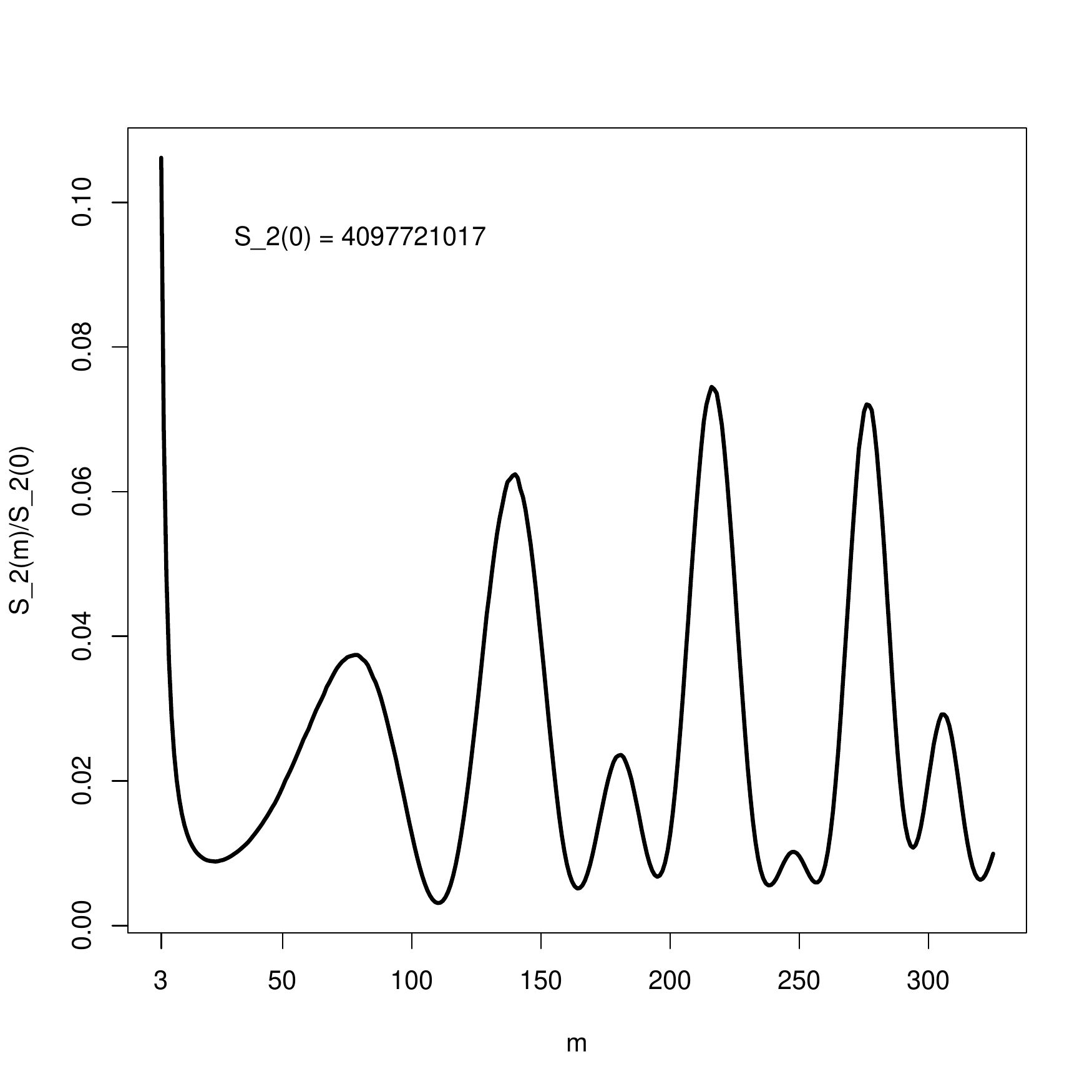}
\end{figure}

\noindent
We computed this shifted moment function for various choices of $m$, $T$, and $H$.
(We also considered similar sums with exponents other than 4, but for simplicity
do not discuss them here.)
Figure~\ref{corrder23} presents some of our results
near the $10^{16}$-th and $10^{23}$-rd zeros, and with $H$ spanning about $10^7$ zeros 
in both cases. 
The figure shows that correlations do exist and persist over long ranges.  
 Also, the shape of $S_2(T,H,m)$ near the $10^{16}$-th zero is similar to that near the $10^{23}$-rd zero, 
 except the former has higher peaks, and covers the range $3\le m\le 222$, as opposed to $3\le m\le 325$, 
 which suggests oscillations scale as $1/\log (T/2\pi)$. 

\begin{figure}[ht]
\footnotesize
\caption{\footnotesize Plot of $M(T,H,\alpha)/M(T,H,0)$, with $H\approx 6.5\times 10^5$, near the $10^{23}$-rd zero, drawn for $\alpha$ a multiple of $0.5$. The dashed line is a sine kernel.}\label{smg2}
\includegraphics[scale=.5]{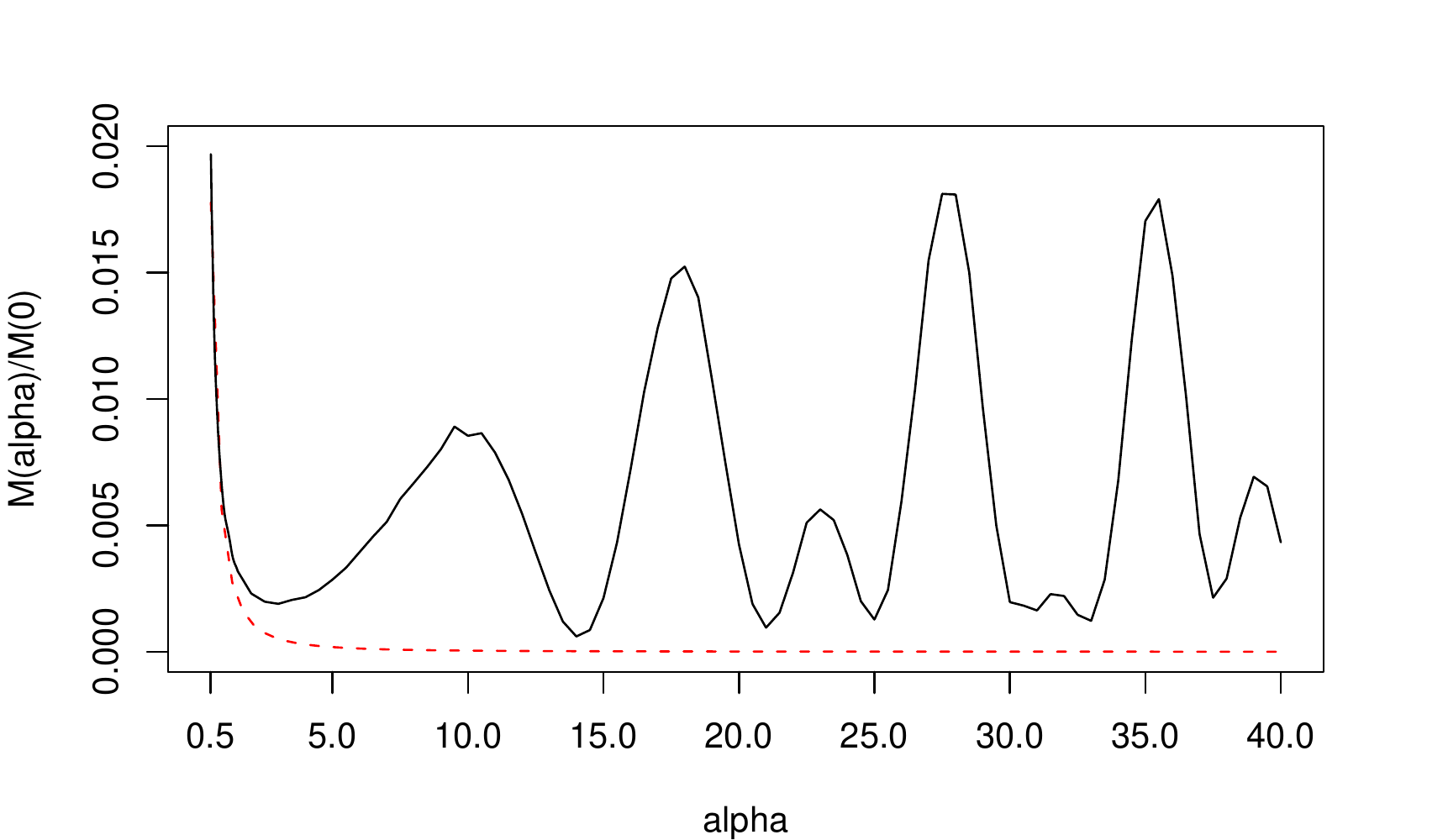}
\end{figure}

We remark the plot of $S_2(T,H,m)$ in Figure~\ref{corrder23} (right plot) 
 is similar to a plot in \cite{HO} of the shifted fourth moment of the zeta function on the
critical line:  

\begin{equation}\label{eq:malpha}
M(T, H; \alpha) := \int_T^{T+H} |\zeta(1/2+it)|^2\,|\zeta(1/2+it +i\alpha)|^2\,dt\,,
\end{equation}
 
\noindent
which we reproduce here in Figure~\ref{smg2} for the convenience of the reader.

% It is not clear why the two figures would be governed by similar asymptotic oscillations.
% \footnote{The initial segments of $S_2(T,H,m)$ and $M(T,H,\alpha)$ look different
%  since they have opposite convexities.}

%> pdf('16moment_shift4_2.pdf')
%> xfit<-seq(3,325,by=1); plot(xfit,cv16_2$acf[xfit+1]/cv16_2$acf[1], type='l', lwd=2.5, xlab='m',ylab='S_2(m)/S_2(0)',xaxt='n'); axis(side=1,at=c(3,seq(50,300,length=6))); text(80,0.15,'S_2(0) = 251745151')
%> dev.off()

\begin{figure}[htbp]
\footnotesize
\caption{\footnotesize Plots of $f(T,H,x)$, defined in (\ref{eq:specf}), using $10^6$ zeros near the $10^{16}$-th zero (upper left), $10^{20}$-rd zero (upper right), and $10^{23}$-rd zero (lower left). The lower right plot is another plot near the $10^{23}$-rd zero, except it uses a different set of  $2\times 10^6$ zeros.}\label{fftder23}
\includegraphics[scale=.3]{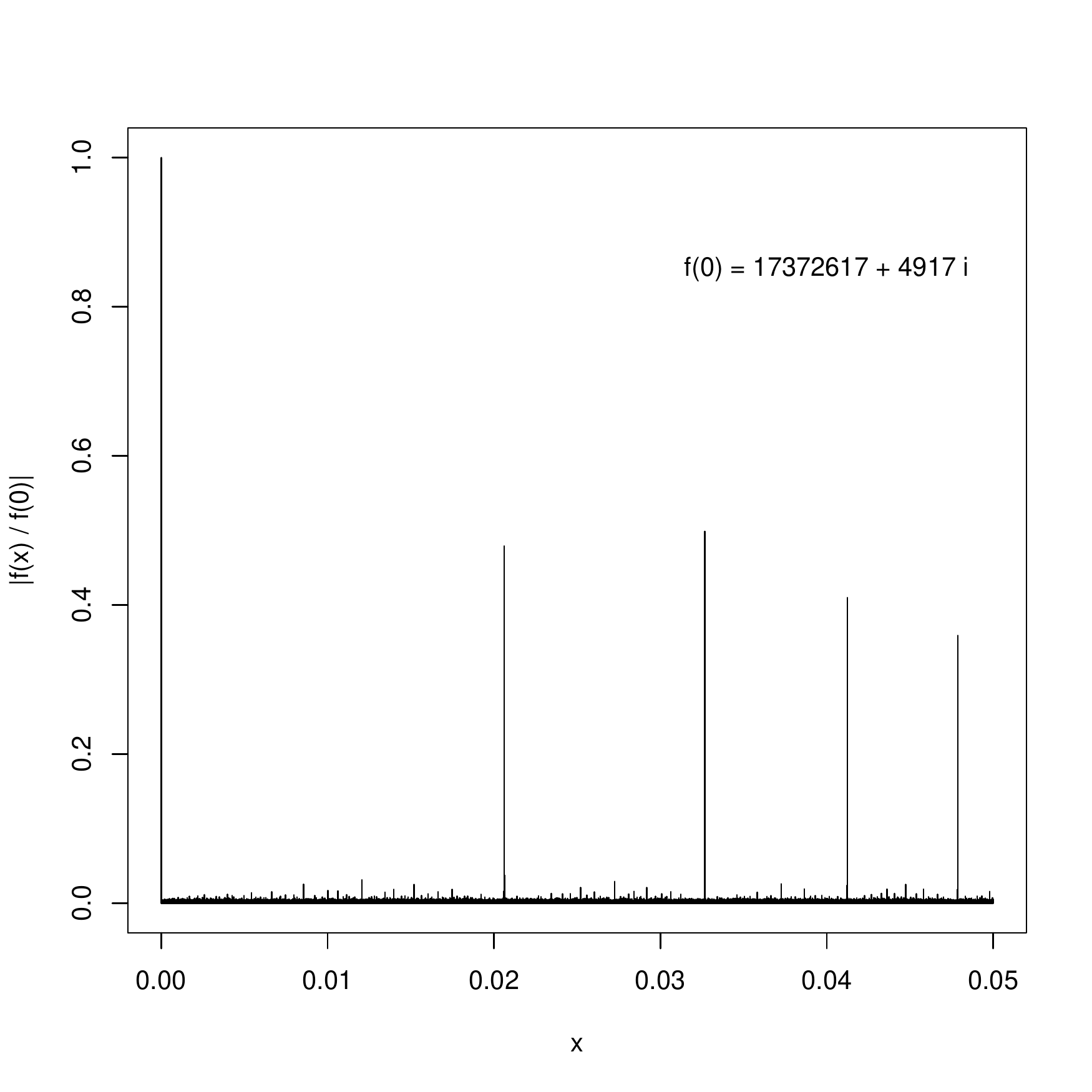}
\includegraphics[scale=.3]{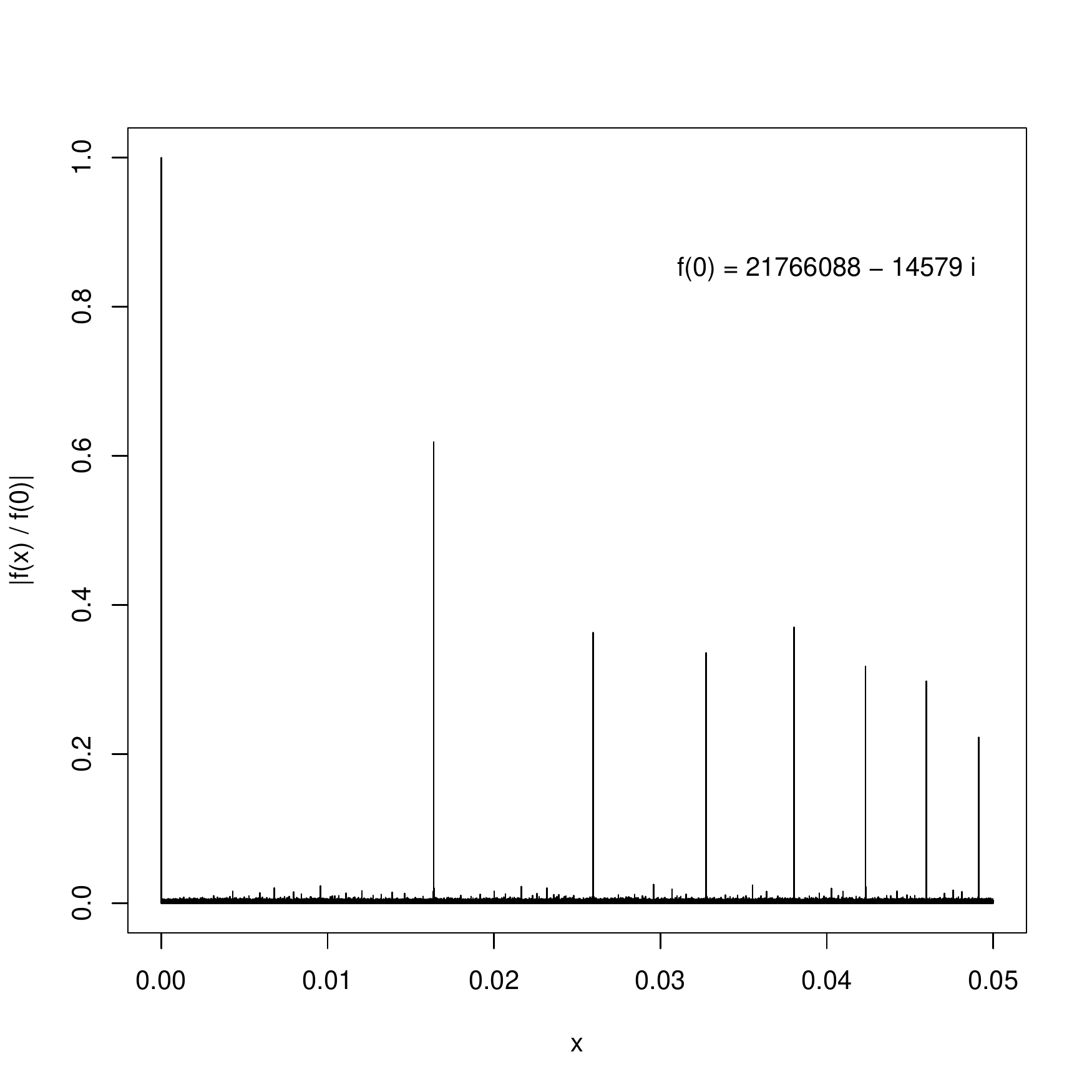}

\includegraphics[scale=.3]{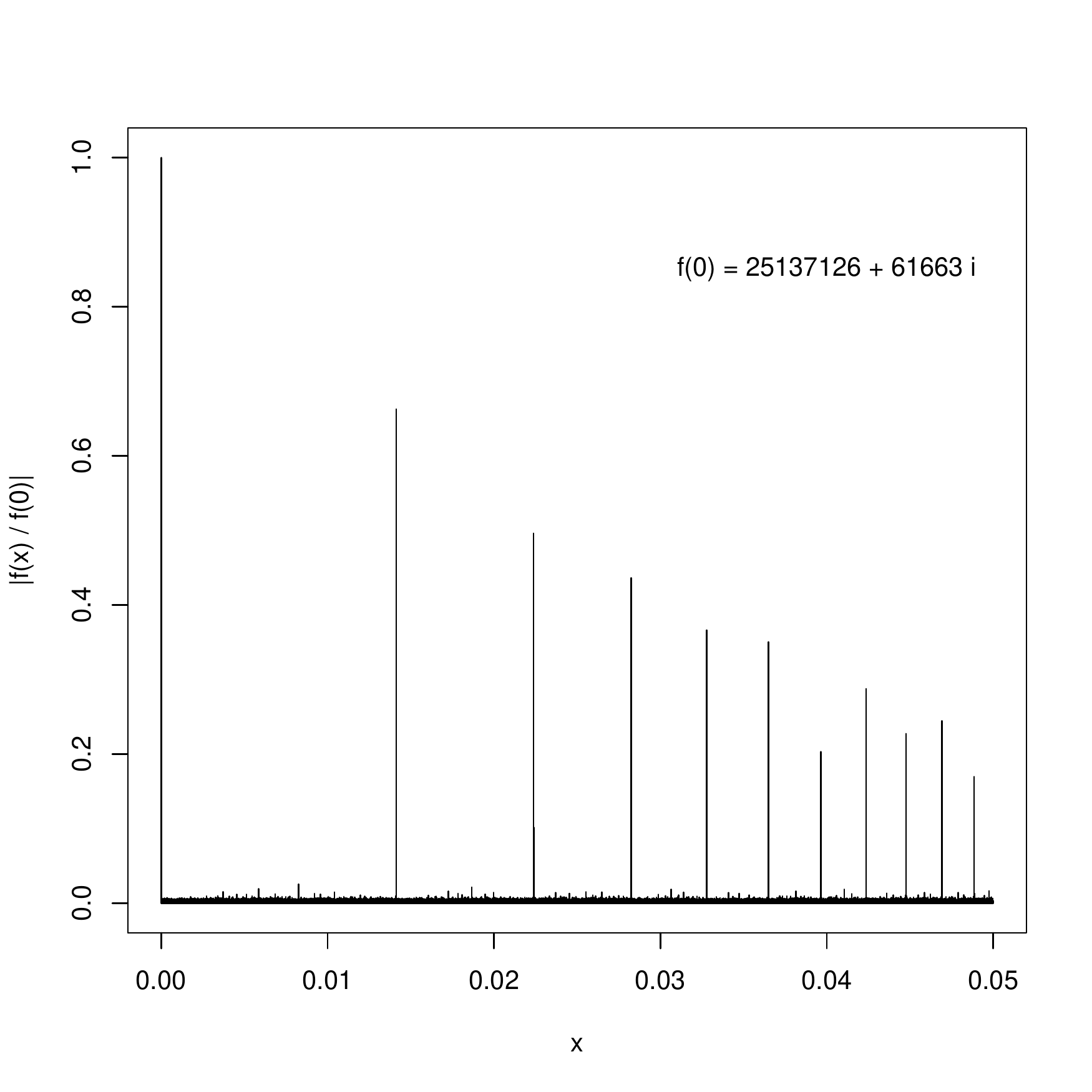}
\includegraphics[scale=.3]{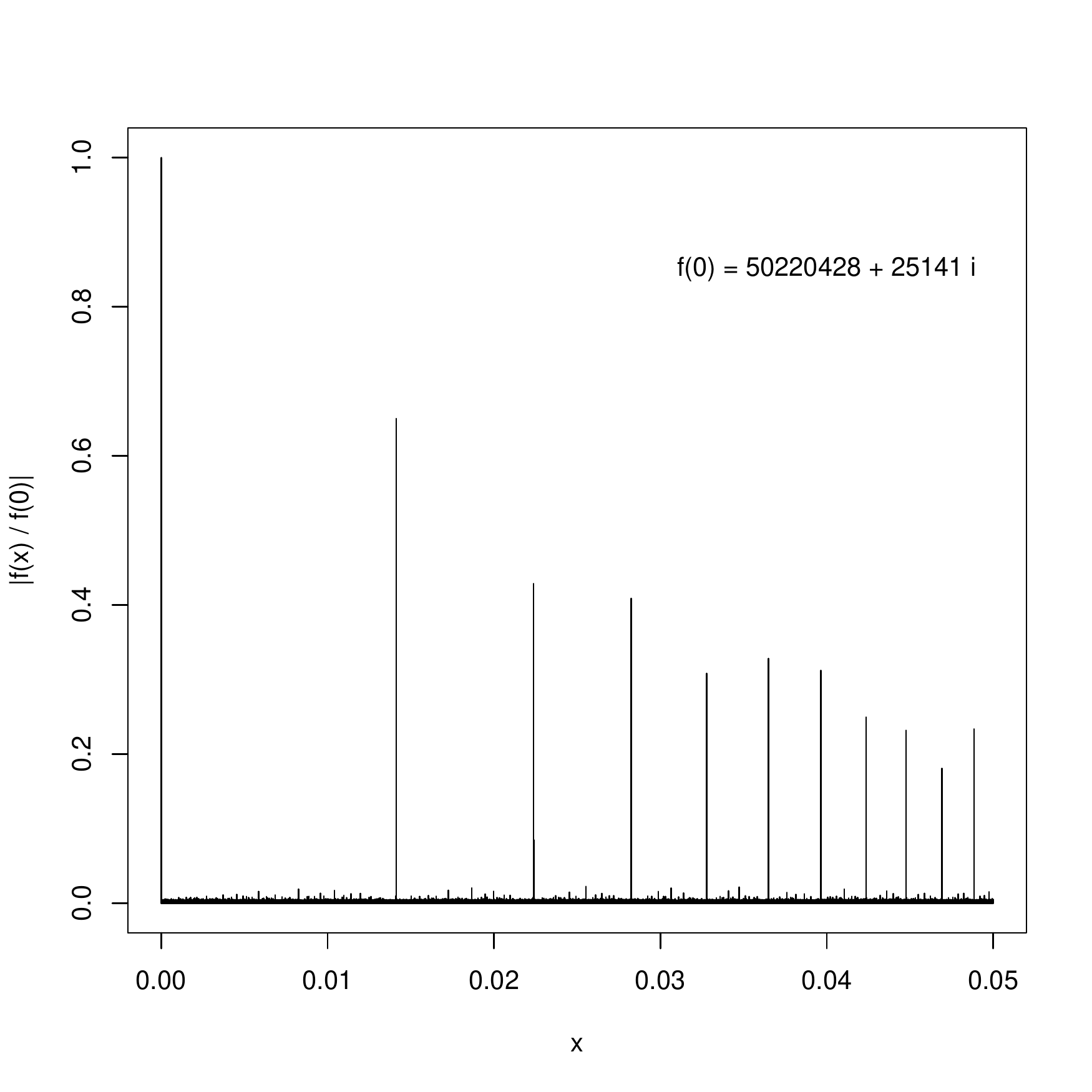}
\end{figure}

To explain observed correlations, we numerically calculated the function:

\begin{equation}\label{eq:specf}
f(T,H,x)=\left|\sum_{T\le \gamma_n\le T+H} \zeta'(1/2+i\gamma_n) e^{2\pi i nx}\right|\,,
\end{equation}

\noindent
which is related to long-range periodicities in $\zeta'(1/2+i\gamma_n)$.
 Assuming the RH, Fujii~\cite{Fu} supplied the
 following asymptotic formula in the case $x=0$: 
  
\begin{equation}\label{eq:fuj}
\sum_{0<\gamma_n\le T} \zeta'(1/2+i\gamma_n)= \frac{T}{4\pi} \log^2\frac{T}{2\pi}+(c_0-1)\frac{T}{2\pi}\log\frac{T}{2\pi} -(c_1+c_0)\frac{T}{2\pi}+O\left(T^{1/2}\log^{7/2} T\right)\,, 
\end{equation}

%\noindent
%where $c_0$ and $c_1$ are the coefficients in the expansion 

%\begin{equation}
%\zeta(s)=(s-1)^{-1}+c_0-c_1(s-1)+\frac{c_2}{2!}(s-1)^2-\frac{c_3}{3!}(s-1)^3+\cdots\,. 
%\end{equation}

\noindent
where $c_0=0.5772\ldots$ (the Euler constant) and $c_1=-0.0728\ldots$.
 Empirical values of $f(T,H,0)$ agree well with formula (\ref{eq:fuj}).
 For example,  with $H$ spanning $10^6$ zeros, we obtain $f(T,H,0)= 21766088 - 14579 i$ 
 near the $10^{20}$-zero,  and we obtain $f(T,H,0)= 25137126+ 61663 i$ near the $10^{23}$-rd zero. 
 But as $x$ increases, $f(T,H,x)$ experiences sharp spikes for certain $x$, as shown in Figure~\ref{fftder23},
 which depicts the segment $0\le x\le 0.05$ (in the remaining portion $0.05< x<1$, the spikes get progressively denser). 
  
The sharp spikes in Figure~\ref{fftder23} show the
 existence of long-range periodicities among values of $\zeta'(1/2+i\gamma_n)$.
These spikes, as well as the correlations described above, 
are not unexpected.  They can be demonstrated to follow from the
properties of the zeta function, by estimating proper contour
integrals.  Such methods were used for continuous averages by Ingham \cite{Ingh}
and even others before him, and for discrete averages over zeros by
Gonek~\cite{Go1} and Fujii~\cite{Fu,Fu2}.  The main step involves integration
of ${\zeta'(s)}^2/{\zeta (s)}$, and estimates of such integrals. 

Applying such methods to ${\zeta '(s)}^2 e^{x s \log\frac{T}{2\pi}} /{\zeta (s)}$  
suggests that the function 
 
\begin{displaymath}
\tilde{f}(T,H,x)=\left|\sum_{T\le \gamma_n \le T+H} \zeta'(1/2+ i \gamma_n) e^{2\pi i \tilde{\gamma}_n x}\right|\,, 
\qquad \tilde{\gamma}_n := \frac{\gamma_n}{2\pi}\,\log\frac{T}{2\pi}\,,
\end{displaymath}

\noindent
experiences large spikes at approximately $x=\log(k) /\log(T/(2\pi))$.
 For by a heuristic argument  
 involving the (very) regular spacing of zeros  
 one expects that $\tilde{\gamma}_n$ in the definition of $\tilde{f}(T,H,x)$ can be replaced by $n$ without too much error 
 (see \cite{Od2} for a similar argument in the context of long-range correlations in zero spacings).
 Therefore, $f(T,H,x)$ should behave similarly to $\tilde{f}(T,H,x)$.\footnote{Indeed, 
  the plots in Figure~\ref{fftder23} are almost unchanged if instead of plotting $f(T,H,x)$ we plot
 $\tilde{f}(T,H,x)$.}
 In particular, we expect the $k$-th spike in Figure~\ref{fftder23} to occur
  at approximately $\log(k) /\log(T/(2\pi))$,  
 and that agrees well with the evidence of the graphs.

\section{Numerical methods}
As usual, define the rotated zeta function on the critical line by
\begin{equation}
Z(t)= e^{i\theta}\zeta(1/2+it)\,,\qquad e^{i\theta(t)}= \left(\frac{\Gamma(1/4+it/2)}{\Gamma(1/4-it/2)}\right)^{1/2} \pi^{-it/2}\,. 
\end{equation}
\noindent
The rotation factor $e^{i\theta(t)}$ is chosen so that $Z(t)$ is real. 
 In our numerical experiments, $t < 1.31 \times 10^{22}$. 

 Since $|Z'(\gamma_n)|=|\zeta'(1/2+i\gamma_n)|$,
 it suffices to compute $Z'(\gamma_n)$.
 To do so, we used
 the numerical differentiation formula (Taylor expansion)
\begin{equation}\label{eq:numdiff}
Z'(t) = \frac{Z(t + h) - Z(t - h)}{2 h} + R(t,h)\,,
\end{equation}

\noindent
where the remainder term in (\ref{eq:numdiff}) satisfies
\begin{equation}
|R(t,h)| \le  \max_{t - h \le t_1 \le t + h}  \frac{|Z'''( t_1 )|}{6} h^2\,,
\end{equation}

\noindent
We chose $h = 10^{-5}$, and approximated the derivative by
\begin{equation}\label{eq:numdiff1}
Z'(t) \approx \frac{Z(t + h) - Z(t - h)}{2 h}\,.
\end{equation}

To evaluate $Z(t)$ at individual
points, we used a version of the Odlyzko-Sch\"onhage algorithm~\cite{OS} 
implemented by the second author~\cite{Od1}.
 If  the point-wise evaluations 
 of $Z(t+h)$ and $Z(t-h)$ 
 via this implementation
 are accurate to within 
 $\pm \epsilon$ each, then the approximation (\ref{eq:numdiff1})
  is accurate to within $\pm (10^5\epsilon + |R(t,h)|)$.
 Numerical tests suggested $\epsilon$ 
 is normally distributed with mean zero and standard deviation $10^{-9}$.
 Therefore, $\epsilon$ is typically around $10^{-9}$.
 Also, varying the choice of $h$ in (\ref{eq:numdiff1}) suggested 
 the approximation is accurate to about 4
 decimal digits with $h=10^{-5}$ and $t\approx 10^{22}$.

In principle, our computations of $\zeta'(1/2+i\gamma_n)$
can be made completely rigorous
by carrying them out in sufficient precision.  
 If one plans on calculating $\zeta'(1/2+i\gamma_n)$ with very high
  precision, however, it will likely be better to first derive a
Riemann-Siegel type formula for $Z'(t)$ itself, with explicit estimates for
 the remainder. Such a formula will be useful on its own as it can be
 be used to check other conjectures about $\zeta'(1/2+it)$.

\section{Conclusions}
Numerical data from high zeros of the zeta function generally
agrees well with the asymptotic results that have been proved,
as well as with several conjectures.  There are some systematic
differences between observed and expected distributions, but
the discrepancies decline with growing heights.

The results of this paper provide additional evidence for the
speed of convergence of the zeta function to its asymptotic
limits.  They also demonstrate the importance of outliers,
and thus the need to collect extensive data in order
to obtain valid statistical results.  The long-range correlations
that have been found among values of the derivative of the zeta
function at zeros can be explained by known analytic techniques.

\end{document}